\newif
\ifconfver
\confvertrue        %declaring conference version true

\ifconfver
    \documentclass[10pt,journal,final,twoside]{IEEEtran}
\else
    \documentclass[11pt,draftcls,onecolumn]{IEEEtran}
\fi
\usepackage{xspace,empheq,fancybox,amsmath,bbm,amssymb,epsfig,subfigure,syntonly,times,amsthm} 
\usepackage{psfrag,color,bm,array,cite}
\usepackage{url,cite,footnote,xspace,syntonly,algorithm,algpseudocode}
\usepackage{verbatim,multirow,slashbox,enumitem}
\usepackage[T1]{fontenc}
\usepackage[top=0.57in, bottom=0.57in, left=0.57in, right=0.57in]{geometry}
\input{mysymbol.sty}

\newtheorem{remark}{\bfseries Remark}

\newcommand{\Lag}{{\mathcal{L}}}
\newcommand{\R}{{\mathbb{R}}}

\newcommand{\bq}{{\mathbf{q}}}

\newcommand{\df}{{\nabla f}}
\newcommand{\dg}{{\nabla g}}

\newcommand{\T}{\top}

\newcommand{\ind}{\mathcal{I}}

\newcommand{\dind}{\widetilde{\nabla}\ind}
\newcommand{\Proj}[1]{\mathcal{P}\{#1\}}
\newcommand{\sigmax}[1]{\sigma_{\max}\{#1\}}

\DeclareMathOperator*{\argmin}{arg\,min}

\newcommand{\tbv}{\widetilde{\bbv}}
\newcommand{\tv}{\widetilde{v}}

\newcommand{\ml}{\color{black}{}}
\def\ws#1{{\color{black}#1}}
\newcommand{\seta}{\widetilde{\eta}}
\newcommand{\sL}{\widetilde{L}}

\begin{document}
%
%%%%%%%%%%%%%%%%%%%%%%%%%%%%%%%%%%%%%%%%%%%%%%%%%%%%%%%%%%%%%%%%%%%%%%
%                                                                    %
%               Paper Title                                          %
%                                                                    %
%%%%%%%%%%%%%%%%%%%%%%%%%%%%%%%%%%%%%%%%%%%%%%%%%%%%%%%%%%%%%%%%%%%%%%

%\title{Communication-Cognizant Hybrid Voltage Control in Power Distribution Networks}
\title{Hybrid Voltage Control in Distribution Networks Under Limited Communication Rates}
\author{\IEEEauthorblockN{Hao Jan Liu},~\IEEEmembership{Student Member, IEEE}, \IEEEauthorblockN{Wei Shi},~\IEEEmembership{Member, IEEE}, and \IEEEauthorblockN{Hao Zhu},~\IEEEmembership{Member, IEEE}
\thanks{\protect\rule{0pt}{0mm}
This material is based upon work supported by %the Department of Energy under Award Number DE-OE000078 and 
the National Science Foundation under Award Number EPCN-1610732 and the Power Systems Research Engineering Center (PSERC). H. J. Liu and H. Zhu are with the Department of Electrical \& Computer Engineering, University of Illinois, 306 N. Wright Street, Urbana, IL, 61801, USA; Emails: {\{haoliu6,haozhu\}@illinois.edu}. W. Shi is with the Department of Electrical \& Computer Engineering, Boston University, 8 Saint Mary's Street, Boston, MA 02215, USA; Email: wilburs@bu.edu.}
}

% % The paper headers
%\markboth{IEEE TRANSACTIONS ON SMART GRID}
%{Liu \MakeLowercase{\textit{et al.}}: Hybrid Voltage Control in Distribution Networks Under Limited Communication Rates}
%\renewcommand{\thepage}{}
\maketitle
\pagenumbering{arabic}

%%%%%%%%%%%%%%%%%%%%%%%%%%%%%%%%%%%%%%%%%%%%%%%%%%%%%%%%%%%%%%%%%%%%%%
%                                                                    %
%                   Abstract                                         %
%                                                                    %
%%%%%%%%%%%%%%%%%%%%%%%%%%%%%%%%%%%%%%%%%%%%%%%%%%%%%%%%%%%%%%%%%%%%%%%%
{\ml \begin{abstract}
		
		Voltage regulation in distribution networks is challenged by increasing penetration of distributed energy resources (DERs). Thanks to advancement in power electronics, these DERs can be leveraged to regulate the grid voltage by quickly changing the reactive power outputs. This paper develops a hybrid voltage control (HVC) strategy that can seamlessly integrate both local and distributed designs to coordinate the network-wide reactive power resources from DERs. \ws{By minimizing a special voltage mismatch objective, we achieve the proposed HVC architecture using partial primal-dual (PPD) gradient updates that allow for a distributed and online implementation}. The proposed HVC design improves over existing distributed approaches  by integrating with local voltage feedback. As a result, it can dynamically adapt to varying system operating conditions while being fully cognizant to the instantaneous availability of communication links. Under the worst-case scenario of a total link outage, the proposed design naturally boils down to a surrogate local control implementation. Numerical tests on realistic feeder cases have been to corroborate our analytical results and demonstrate the algorithmic performance.
		
		%The operational paradigm of power distribution networks has witnessed significant transformations with increasing penetration of distributed energy resources (DERs). Voltage regulation, an important distribution operational task, has been greatly challenged by the intermittency and variability of DERs. Thanks to advancement in power electronics, these DERs can be leveraged to regulate the grid voltage by quickly changing the reactive power outputs. This paper develops a hybrid voltage control (HVC) strategy that can seamlessly integrate both local and distributed control designs to optimally coordinate the network reactive power of DERs. By designing a special voltage mismatch objective, we achieve the proposed HVC architecture using the partial primal-dual gradient (PPD) algorithm that can allow for decentralized online implementations. The resultant HVC design improves over existing distributed methods by integrating with local voltage feedback. It can dynamically adapt to varying system operating conditions while being fully cognizant to the instantaneous availability of communication links. Under the worst-case scenarios of a total link outage, the proposed design naturally boils down to a surrogate local control implementations. Numerical tests on realistic feeder cases are presented to corroborate our analytical results and demonstrate the algorithmic performance.
\end{abstract}}

%\begin{IEEEkeywords}
%\end{IEEEkeywords}

%\newpage

%%%%%%%%%%%%%%%%%%%%%%%%%%%%%%%%%%%%%%%%%%%%%%%%%%%%%%%%%%%%%%%%%%%%%%
% %%
% %%      Section: Intro %%
% %%%%%%%%%%%%%%%%%%%%%%%%%%%%%%%%%%%%%%%%%%%%%%%%%%%%%%%%%%%%%%%%%%%%
%%%
\section{Introduction}\label{sec:intro}
Recent proliferation of distributed energy resources (DERs) such as solar generation and storage devices  can potentially cause some rapid voltage fluctuations in distribution networks. 
{\ml 
	%These unexpected fluctuations are mainly due to high intermittency and variability of these resources \cite{Car_tps08, Jaha_tps13}. Additionally, overvoltage scenarios can arise under the reverse power flow from distributed generations. Traditional voltage regulation devices are operated on a minute- or even hourly-basis. Thus, the fast system dynamics along with overvoltage issues challenge the operational goal of maintaining the voltage magnitude to be within $\pm 5 \%$ of its nominal value \cite{Tur_11}. 
	A promising technology to tackle this challenge is via advanced inverter control design. The fast-acting power electronic-interfaced DERs can support the distribution system voltage regulation objective by controlling their \textit{reactive power} (VAR) outputs.}

The voltage control problem can be viewed as a special case of the optimal power flow (OPF) one that minimizes the voltage mismatch. {\ml This OPF-based approach requires the availability of full network-wide information at a centralized location \cite{Far_pesgm12}.} %Hence, such architecture generally suffers from the limitations on communication infrastructure and resources in distribution networks. 
It is also possible to develop local control strategies that only use local voltage magnitude information \cite{Jaha_tps13, Tur_proc11}. Due to lack of information exchange, local designs could suffer from instability issues and sub-optimality \cite{hzml_tps15}.
For reduced communication complexity yet globally optimal performance, several distributed optimization based techniques using information exchanges among neighboring buses have been proposed in \cite{Robbins_tps16,Rob_tps13,dall_tsg13,lanl,hlwshz_tsg17}; see \cite{Hatz_TSG17} for a review of recent distributed and decentralized methods. Albeit a distributed control framework can more effectively coordinate network-wide VAR resources, its performance highly depends on communication capabilities of distribution networks, in particular the rates of information exchange as pointed out by \cite{Hatz_TSG17}. While the effects of bandwidth limits and quantized messages have been investigated in \cite{Li_arxiv17} for distributed voltage control, more focus has been on the analysis of low-rate communications, or equivalently, asynchronous control updates. For example, \cite{Dall_cdc16} has designed an asynchronous decentralized algorithm where each DER controller flexibly incorporates the fast incoming local information with the low-rate control signal sent by a centralized aggregator. It has been shown in \cite{hlwshz_tsg17} that distributed voltage control updates can be also executed in an asynchronous fashion.  Nonetheless, existing distributed control approaches would fail to work for the worst-case scenario of a total communication outage, under which they need to freeze all updates and cannot incorporate any local voltage information.  As distribution systems {\ml continue to witness} low-rate communication networks, it is imperative to integrate both local and distributed control frameworks to achieve the dual objectives in terms of adaptiveness to various communication rates and globally optimal voltage regulation performance. %, facilitating a sustainable development of the future smart grid. 

The present paper develops a hybrid voltage control (HVC) strategy that can dynamically adapt to varying system operating conditions while being fully cognizant to the instantaneous availability of communication links. To cope with practical communication limitations, the proposed HVC scheme consists of both distributed and local control architectures and does not require a centralized authority. We formulate the network-wide VAR optimization problem by linearizing the power flow model for analysis purposes only. The resultant quadratic programming problem is solved by a \textit{partial} primal-dual gradient (PPD) algorithm in discrete-time domain, which is a variant of the basic primal-dual (sub)gradient method, see e.g., \cite{Nedic_jota09}. We provide the analysis in step-size choices that can guarantee convergence. We further use the PPD-based solver to design online HVC strategy where each bus can integrate both the local voltage measurement and the communication information shared by neighboring buses. Although a linearized model has been adopted to the algorithmic development and analysis, performance of the proposed HVC design has been verified using the full ac power flow model {\ml for unbalanced and lossy distribution networks.} 

Compared to existing approaches on voltage control design, the main contributions of our HVC design are three-fold. First, it explicitly accounts for the VAR limits by using the projection operator. Due to the discontinuity of projection mapping, general Krasovskii's methods \cite{Feijer_cdc09} for analyzing the stability of primal-dual gradient flow method would not hold. To tackle this problem, we have expressed the operation as a subgradient step featured by an indicator function in order to establish the stability of the PPD-based HVC design in Sec. \ref{sub:con_a}. Second, the HVC design only requires each bus to {\ml measure} its local voltage magnitude and communicate it to neighboring buses. The sensing requirement and communication overhead are minimal compared to most (de)centralized strategies. Last but not least, our hybrid design can integrate both the neighboring bus voltage information and local voltage measurements regardless of the communication link conditions.  %{\ml Albeit this hybrid design may potentially increase the implementation complexity, }it enjoys the benefits that can be offered by both local and distributed control designs. 
{\ml This way, the HVC design is cognizant to the instantaneous availability of communication links while effectively tracking the globally optimal VAR setting.} Interestingly, although the HVC updates have been developed using the distributed PPD-based solver, it would boil down to a surrogate local voltage control update during a total communication outage. Under this worst-case scenario, satisfactory performance can still be achieved as it responds to local voltage variation.

This paper is organized as follows. Section \ref{sec:model} presents the power flow model for distribution networks and formulates the specially designed voltage optimization problem. The PPD-based HVC algorithm is developed in Section \ref{sec:method}, along with its online implementation and communication-cognizant design. Numerical test results using realistic feeders and real-time-series data are offered in Section \ref{sec:num}. Conclusions are presented in Section \ref{sec:con}.

%%%%%%%%%%%%%%%%%%%%%%%%%%%%%%%%%%%%%%%%%%%%%%%%%%%%%%%%%%%%%%%%%%%%%%
% %%
% %%      Section: System Modeling and Problem Statement%%
% %%%%%%%%%%%%%%%%%%%%%%%%%%%%%%%%%%%%%%%%%%%%%%%%%%%%%%%%%%%%%%%%%%%%
%%%
\section{System Modeling and Problem Statement} \label{sec:model}
%
%\begin{figure}[tb]
%	\centering
%	\vspace{-3pt}
%	\includegraphics[width=1\linewidth,clip = true]{radial}
%	\caption{A radial distribution network with bus and line associated variables.}
%	\label{fig:radial}
%	\vspace{-4mm}
%\end{figure}
%
We consider a  distribution network given by a graph with the set of buses $\ccalN:=\{0,...,N\}$ and the set of line segments $\ccalE:=\{(i,j)\}$.
%; see Fig. \ref{fig:radial} for a radial network configuration. 
Per bus $j$, let $v_{j}$ denote its voltage magnitude, and $p_j$ $(q_j)$ represent the active (reactive) power injection, respectively. All network quantities are in per unit (p.u.). A constant reference bus voltage $v_0$ is assumed for the point of common coupling. For each line $(i,j)$, we denote $r_{ij}$ and $x_{ij}$ as its resistance and reactance in addition to $P_{ij}$ and $Q_{ij}$ as the active and reactive power flow from $i$ to $j$, respectively. 
The so-termed LinDistFlow model has been developed in \cite{Baran_tpd89} to linearize the power flow model, assuming negligible line losses and almost flat voltage. Its accuracy can be numerically corroborated by several recent work \cite{Far_cdc13, Sulc_tec14,hzml_tps15,hlwshz_tsg17}. Per bus $j$, we consider the controllable VAR $q^\mathrm{g}_j:=q_j+q_j^\mathrm{c}$ as our control input where $q_j^\mathrm{c}$ is the VAR  consumption. It is shown that one can relate the voltage $v_j$ to the controllable VAR $q_j^\mathrm{g}$ as the following \cite{hzml_tps15}:
\begin{equation}\label{ldfVm2}
\begin{array}{c}
\sum_{i \in \mathcal{N}_{j}} B_{ji} v_i = q_j^\mathrm{g} + w_j,~\forall j \in \ccalN 
\end{array}
\end{equation}
where $\ccalN_j:=\{i | (i,j) \in \ccalE \} \cup \{j\} \subseteq \ccalN$ contains bus $j$ and all of its neighboring buses; the quantity $w_j$ captures the system operating conditions, i.e., the effect of $p_j$ on $v_j$ when $q_j^\mathrm{g}=0$. By concatenating all scalar variables into vectors and replacing $\bbq^\mathrm{g}$ by $\bbq$ for notational convenience, 
%To better present our paper, hereafter we will frequently use bold letters to denote the concatenation (vector form) of corresponding scalar variables. In addition, we replace $\bbq^\mathrm{g}$ by $\bbq$ for notational convenience.
\eqref{ldfVm2} can be represented in a compact form $\bbB \bbv = \bbq + \bbw$. Matrix $\bbB$ is the \textit{Bbus matrix} used in the dc power flow model (see e.g., \cite[Sec. 6.16]{Wood_book}).  We reserve $\bbX=\bbB^{-1}$ for our future use. By definition, matrix $\bbB$ is a reduced, weighted graph Laplacian matrix (full rank) and thus has a unique sparsity pattern based on the network topology. Similarly, a linearized model using the graph-based matrices also holds for multi-phase unbalanced networks \cite{Gan_pscc14}. For simplicity, we will present the proposed algorithm using \eqref{ldfVm2}.

We aim to optimize the network $\bbq$ contributed from local inverters to {\ml achieve a desired voltage profile $\bbmu$ such that  $\bbv \to \bbmu$. This is similar to the goal of secondary voltage control design in transmission and microgrid networks to effectively coordinate the network VAR resources and enhance the voltage stability \cite{su_tsg16,lu_tetcas15}. Such preferred profile $\bbmu$ can be adjusted depending on certain operational specifications for the network, e.g., conservation voltage reduction implementations would be most effective if the voltage profile is flat ($\bbmu=\mathbf{1}$) \cite{Liu_sgcomm14}.}
We first formulate the problem under a static setup where the operating condition vector $\bbw$ is constant to develop the HVC scheme. {\ml As detailed soon in Sec. \ref{sec:method}-B, we will extend it to an online design by dynamically updating $\bbw$ using neighboring voltage magnitude measurements.} To this end, the voltage control objective of minimizing $f_1(\bbv) = \frac{1}{2} \|\bbv-\bbmu\|^2$ {\ml is introduced to} improve the  system voltage level by coordinating network-wide VAR resources and hence providing the \textit{globally} optimal VAR setting. Minimizing $f_1$ subject to \eqref{ldfVm2} could be solved collaboratively by all VAR resources through communication between buses using distributed optimization techniques \cite{Rob_tps13,dall_tsg13,lanl}. Nonetheless, these distributed control schemes would fail to work under a total communication outage. To make the control design more robust, we introduce another \textit{weighted} objective function $f_2(\bbv)=  \frac{1}{2} \|\bbv-\bbmu\|_{\bbB}^2$ where we define the weighted norm $\|\bby\|_\bbB ^2 := \bby^\T\bbB\bby$ for any vector $\bby$. Interestingly, it turns out that minimizing $f_2$ subject to \eqref{ldfVm2} can be tackled through a totally local control architecture, i.e., voltage droop control scheme in \cite{hzml_tps15}. {\ml Albeit this weighted objective inevitably results in a sub-optimal VAR setting under limited VAR resources, it turns out this objective can allow for communication-free updates using only local voltage measurements. 
	%its communication-free feature is rather attractive considering practical cyber resource constraints. 
	We form the total objective using both $f_1$ and $f_2$, such that the resultant design would enjoy the dual benefits of being globally optimal and communication-cognizant.}
%we integrate both $f_1(\bbv)$ and $f_2(\bbv)$ into the whole objective, 
%expecting that the yielded hybrid strategy is more effective than the pure local design. Though compromising the optimal design a bit, the new hybrid strategy would be communication-cognizant and robust to the worst-case scenario of a total communication link outage.} 
%Furthermore, compared to the traditional paradigm of maintaining the voltage within limits, these aforementioned objectives of minimizing voltage error mismatch can improve the network-wide voltage profile by coordinating VAR resources across the network. To allow some tolerance in the voltage mismatch error, it is possible to use other convex error penalty functions such as the Huber loss \cite[Ch. 7]{huber_robust81} instead of the squared $\ell_2$ norm in $f_1$ and $f_2$. 
Combining the attractive features of both distributed and local control frameworks, we \ws{cast} the HVC problem as 
\begin{subequations}
	\label{intvar}
	\begin{align}
	\{\bbv^\star,\bbq^\star\}\!:=\! \argmin_{\bbv,\bbq} &~ h(\bbv,\bbq) := f_1(\bbv) + \gamma f_2(s(\bbq)) \label{wvobj} \\
	\mathrm{subject~to}&~\bbB \bbv = \bbq + \bbw \label{bvq} \\
	&~~\underline{\bbq}\leq\bbq\leq\overline{\bbq} \label{invlim}
	\end{align}
\end{subequations}
where the function $s(\bbq):=\bbX(\bbq+\bbw)=\bbv$ based on the physical system couplings\footnote{Note that here we assume each bus $j$ has the knowledge of its own operating conditions $w_j$.}. {\ml Hence, both $\bbv$ and $\bbq$ are the decision variables in the problem \eqref{intvar}.} As detailed soon in Sec. \ref{sec:method}, inclusion of $\bbq$ into the formulation would lead to a special feature where the local voltage measurement is {\ml part of the instantaneous gradient direction with respect to $\bbq$, which  allows for separable problem structure and provides the VAR control law.} Hence, we can deal with this part of objective without communication. Note that problem \eqref{intvar} is still convex because $\bbB$ is positive definite \cite{hzml_tps15}. 
%In addition, \eqref{bvq} is equivalent to the power flow constraints in \eqref{ldfVm2} due to the sparse structure of $\bbB$.  
As $\bbB$ is a weighted graph Laplacian matrix, it has a unique structure according to the network topology, i.e., $B_{ij}=B_{ji}=0, \forall (i,j) \notin \mathcal{E}$. Accordingly, coupling involving only the  neighborhood $\ccalN_j$ in \eqref{bvq} is instrumental for the HVC design in Sec. \ref{sec:method}. We also introduce the importance factor $\gamma > 0$ to account for the importance of the local design. {\ml To enhance the performance under the extreme total communication outage scenario, $\gamma$ should be chosen as large as possible to facilitate the communication-free feature of $f_2$. Meanwhile, if the communication link quality is very high, one can decrease the $\gamma$ value toward $0$ to better achieve the unweighted voltage mismatch minimization criterion.}
%This way, a trade-off between control performance under a total communication outage and voltage mismatch level always exists.} 
Last, \eqref{invlim} is the aggregation of box constraints $\underline{q}_j \leq q_j \leq \overline{q}_j,~\forall j \in \ccalN$. This restriction either comes from the inverter apparent power limit or depends on certain inverter power factor limit. Considering the objective in \eqref{wvobj}, one can see that the voltage mismatch error is a trade-off between the distributed (communication-involved) and local (communication-free) control objectives. Suppose that every bus has unlimited VAR capability (i.e., no presence of \eqref{invlim}), 
%{\ml we rewrite the objective function (5a) based on $\bbv = s(\bbq):=\bbX(\bbq+\bbw)$, as given by $$h(\bbq):=\frac{1}{2} \|\bbX(\bbq+\bbw)-\bbmu\|^2+\frac{\gamma}{2} \|\bbX(\bbq+\bbw)-\bbmu\|_{\bbB}^2.$$ Taking the gradient direction with respect to $\bbq$ and setting it to zero, we have $(\bbX+\gamma\bbI)[\bbX(\bbq+\bbw)-\bbmu]=\bb0$. Clearly, the minimizer without the VAR limit constraints $\bbq^\star=\bbB\bbmu-\bbw$ corresponding to the same minimizer of $f_1$.} 
\eqref{intvar} with arbitrary positive (including infinity) $\gamma$ would give the same optimal solution of $f_1$. This implies that the optimal solution of \eqref{intvar} has the potential to closely approximate the globally optimal VAR setting obtained by minimizing $f_1$ under abundant VAR resources.

\section{Hybrid Voltage Control} \label{sec:method}
This section presents our proposed hybrid voltage control (HVC) framework. \ws{Combining the distributed and local control features, we aim to solve \eqref{intvar} by adopting a typed primal-dual gradient  algorithm. Many variants of such algorithm have been studied in both continuous-time domain (see e.g., \cite{Feijer_auto10,Li_tcns14}) and discrete-time domain (see e.g., \cite{Nedic_jota09,zhang2008impact}).} Based on the PPD, we design the HVC scheme that needs only to measure and incorporate the dynamic voltage magnitude. As detailed soon, the proposed feedback approach is very different from existing distributed control schemes in power systems since most of them are developed as a static optimization problem and overlook the communication imperfectness and online implementations. Most importantly, our HVC scheme would boil down to a surrogate local voltage control problem under the worst-case scenario of a total communication failure. Thus, it enjoys a satisfactory performance by having the cognizance to varying communication scenarios.

Introducing the Lagrangian multiplier $\bblambda$ for the equality constraints in \eqref{bvq}, we obtain the following Lagrangian function for the given static optimization problem \eqref{intvar}:
\begin{equation}\label{lagran}
\begin{array}{rcl}
\underset{\bbq \in \ccalQ} {\ccalL(\bbv,\bbq,\bblambda)} 
& = & \frac{1}{2} \|\bbv - \bbmu\|^2 + \frac{\gamma}{2} \|\bbX (\bbq+\bbw)-\bbmu\|_\bbB^2\\
&   & + \langle\bblambda,\bbB\bbv - \bbq - \bbw\rangle
\end{array}
\end{equation}
where we define $\ccalQ := \{\bbq\big|\bbq \in [\underline{\bbq},\overline{\bbq}] \}$ and $\langle \cdot,\cdot \rangle$ represents the inner product. \ws{The goal becomes seeking the saddle point of \eqref{lagran} in a physically implementable way. A popular and efficient method to numerically find the saddle point is to instead construct an augmented Lagrangian and perform alternating minimization over primal variables followed by a (sub)gradient ascent update over dual variables (see e.g., \cite{boyd2011distributed} for the introduction of the ADMM algorithm). Unfortunately, such algorithm would either involve a centralized computing or require multi-hop communication graph. Accordingly, this implementation does not work for the problem of interest here. Hence, we instead focus on the Lagrangian function \eqref{lagran} and adopt the \textit{partial} primal-dual (PPD) algorithm, which admits to a sparse communication network and turns out to cope up with total communication link failures. The "partial" update property will become clear soon, with the algorithmic design detailed here.}

\ws{Using the superscript to denote the iteration number, we initialize $\bblambda^0 = \bb0$ and $\bbq^0$ to be the latest VAR output setting. The PPD updates at the $(k+1)$-st iteration per bus $j$ consists of the following three steps:}
\begin{enumerate}[label=\bfseries (S\arabic*)]
	\item 
	Update $\bbv$: For given $\bbq^k \textrm{ and }\bblambda^k$, $\bbv$ is updated by solving
	\begin{equation}\label{v-update}
	\bbv^{k+1}=\argmin_{\bbv} \ccalL(\bbv,\bbq^k,\bblambda^k)
	\end{equation}
	which is an unconstrained convex quadratic programming. Thus, a closed-form solution exists. Since the objective in \eqref{v-update} is decoupled in $v_j$, the update per bus $j$ is obtained accordingly as
	\begin{equation}
	\label{vbarup}
	v_j^{k+1} = - \sum_{i \in \ccalN_j} B_{ji} \lambda_i^k+\mu_j.
	\end{equation}
	\ws{This is one of the main differences to the classical gradient flow algorithm which would have used (the discretization of) $\dot{\bbv}=-\frac{\partial}{\partial \bbv}\ccalL(\bbv,\bbq,\bblambda)$. Hence, the term ''\textit{partial}'' is adopted here in the name of our algorithm. Similar strategy of ``partial gradient update'' has also been explored in \cite{Li_tcns14} for a continuous-time algorithm whereas we focus on a discrete-time update. Note that the variable $\bbv^k$ estimates the network voltage magnitude based on the power flow model \eqref{ldfVm2}.}
	%Hence, $\bbv^k$ is only adopted to conduct the numerical computations in the ensuing $\bblambda$ update.}
	%whereas $\tdbbv^k$ serves as the network physical voltage magnitude.}
	%Additionally, it is suggested that $\bblambda$ is initialized to zero. As detailed in Sec. \ref{sec:AHVC}, when there is a total communication outage in the network, the PPD-based algorithm would boil down to a pure droop control design based on the current value of $\bblambda^k$ (see \eqref{qjupdate}), i.e., minimizing $f_2$ subject to voltage limits \eqref{invlim}.
	\item 
	Update $\bbq$: Defining step-size $\alpha > 0$, we perform a gradient-projection-based update on $\bbq$, as given by
	\begin{align}
	\label{qupdate}
	\bbq^{k+1} := \Proj{\bbq^k - \alpha \nabla_{\bbq} \ccalL(\bbv^{k+1},\bbq^k,\bblambda^k)}
	\end{align}
	where the projection operator $\Proj{\cdot}$ bounds any input to be within $\ccalQ$. As $\bbq^{k+1}\in [\underline{\bbq},\overline{\bbq}]$ always holds, we use the latest iteration of \eqref{qupdate} to be the network VAR control signal. Under this setting, the gradient direction in \eqref{qupdate} for any pair $(\bbq^k,\bblambda^k)$ becomes
	\begin{align}
	\hspace{-1em}\nabla_{\bbq}\ccalL(\bbv^{k+1},\bbq^k,\bblambda^k) := &~\gamma {\ml \bbX\bbB} [\bbX (\bbq^k+\bbw)-\bbmu)] - \bblambda^k \nonumber \\
	\approx &~\gamma  (\tbv^k - \bbmu) - \bblambda^k \label{qupdate_grad}
	\end{align}

	where $\tbv^k \approx \bbX (\bbq^k+\bbw)$ is the instantaneous network voltage profile at time $k$ following from \eqref{ldfVm2}. Note that $\tilde{v}_j^k$ can be obtain locally by the voltage measurement unit $j$. {\ml 
		With $\bbX \bbB = \bbI$, the gradient function \eqref{qupdate_grad} can be partially computed using the voltage measurement $\tdbbv^k$. By weighting the objective norm of $f_2$ with the positive definite matrix $\bbB$, its gradient direction decouples from the full vector $\bbq$ and can be computed separately with a pair of the voltage magnitude measurement $\tilde{v}_j^k$ and $\lambda_j^k$ at each bus.}
	%As matrix $\bbX$ is a positive definite, its inverse matrix $\bbB$ is unique and positive definite as well. Thus, any other positive definite weighting matrix would not necessarily lead to this decoupling feature.}
	%Due to the weighted norm in $f_2$, the gradient direction completely decouples from the full vector $\bbq^k$ and can be computed separately with each bus pair $(\tilde{v}_j^k,\lambda_j^k)$ as $\bbX \bbB = \bbI$.} 
	Thanks to the separability of box constraints, under arbitrary initialization $q_j^0 \in [\underline{q}_j,\overline{q}_j]$, the update \eqref{qupdate} for each bus $j$ becomes
	\begin{align}
	\label{qjupdate}
	q_j^{k+1} := \ccalP_j\{q_j^k - \alpha \left[ \gamma (\tv_j^k-\mu_j)-\lambda_j^k \right]\}
	\end{align}
	where $\ccalP_j$ denotes the projection at bus $j$ to $[\underline{q}_j,\overline{q}_j]$. 
	{\ml 
		%The advantage of this projection operation for $\bbq$ over the exact minimization of the \eqref{lagran} in (S2) is that it can be implemented locally. 
		Albeit $\ccalP_j$ is constant in this static setup, it would vary according to the inverter limits and instantaneous active power generations under online implementations.} Due to the physical power flow couplings, the local voltage measurement-based gradient direction contains the most up-to-date network-wide information. As detailed soon, this unique feature combined with local dual variable $\lambda_j$ results a robust, communication-cognizant HVC design.
	\item 
	Update $\bblambda$: For a given step-size $\beta$, each multiplier is linearly updated per bus j using the iterative residual of the equality constraints which is equivalent to the gradient $\nabla_{\bblambda}\ccalL(\bbv^{k+1},\bbq^{k+1},\bblambda^k)$. Thus, we have 
	\begin{align}
	\label{lambup}
	\hspace{-2em}\lambda_j^{k+1} := \lambda_j^k + \beta \left[\sum_{i \in \ccalN_j} B_{ji} v_i^{k+1} - q_j^{k+1} - w_j \right].
	\end{align}
\end{enumerate}
%Each iteration consists of the local computations (S1-S3) per bus $j$ in addition to information exchanges among neighboring buses in (S1) and (S3), which leads to a fully distributed HVC design. 
The PPD-based iterations in (S1-S3) constitute the basis for our proposed hybrid voltage control design. 
%
%\begin{remark} \label{rmk:multi} (Multi-phase Networks.) {\rm
%		The PPD-based HVC algorithm can be generalized for the multi-phase network problem. It is possible to obtain the matrix form of \eqref{ldfVm2} for multi-phase networks by further assuming a nearly balance system \cite{Gan_pscc14,Robbins_tps16,kekatos_tps16}. Each scalar variable per bus $j$ for the single phase network is converted to a $3 \times 1$ vector containing all three-phase variables. In addition, the network topology-based Laplacian matrix still holds and thus shares similar sparse structures as $\bbB$, where entries are zero for buses not connected by a line segment. Consequently, this allows for the HVC design in (S1)-(S3) to be extended for unbalance multi-phase networks.}
%\end{remark}

\subsection{Convergence Analysis}
\label{sub:con_a}
Without the projection operator in the $\bbq$-update (in the case $\mathcal{Q}=\R^N$), the convergence analysis of the PPD algorithm would boil down to the stability issue of a discrete-time linear system, whose analysis could be completed in a few lines of words. We would like to treat the system as a linear dynamical system (since linear system is usually considered as easy to handle), allowing us to provide a comprehensive performance analysis. However, with the presence of the projection operator, it becomes a so-called saturated linear system and its stability could be hard to determine \cite{blondel2001stability}. To tackle this problem, we resort to a convex analysis approach. 

The PPD algorithm actually seeks the saddle point of the Lagrangian function \eqref{lagran}. To make our analysis more unified, let us recast the Lagrangian function \eqref{lagran} into a general form
\begin{equation}\label{eq:general}
\Lag(\bbv,\bbq,\bblambda)=f(\bbv)+g(\bbq)+\ind(\bbq)+\langle\bblambda,\bbB\bbv-\bbq-\bbw\rangle
\end{equation}
where $f$ and $g$ are both assumed to be strongly convex and have Lipschitz gradients. In addition, the indicator function
\[
\ind(\bbq)=\left\{
\begin{array}{ll}
0, &\text{if $\bbq\in{\ccalQ}$},\\
+\infty, &\text{if $\bbq\notin{\ccalQ}$}\\
\end{array}
\right.
\]
equivalently accounts for the constraint $\bq\in\mathcal{Q}$. This way, \eqref{lagran} becomes a special case which uses $f(\bbv)\triangleq f_1(\bbv)=\frac{1}{2} \|\bbv - \bbmu\|^2$ and $g(\bbq)\triangleq f_2(s(\bbq))=\frac{\gamma}{2} \|\bbX (\bbq+\bbw)-\bbmu\|_\bbB^2$. By definition, $g(\cdot)$ is a $\eta$-strongly convex function with $L$-Lipschitz gradient, i.e., it holds that
\[
\langle \bbx-\bby,\dg(\bbx)-\dg(\bby)\rangle\geq\eta\|\bbx-\bby\|_2^2,\ \forall \bbx,\bby
\]
with some positive constant $\eta$, and
\[
\|\dg(\bbx)-\dg(\bby)\|_2\leq L\|\bbx-\bby\|_2,\ \forall \bbx,\bby,
\]
with some positive constant $L$. Here, $\langle\cdot,\cdot\rangle$ is the standard inner product. Similarly, $f(\cdot)$ is a $c$-strongly convex function and $\ind(\cdot)$ is a convex function. 

The saddle point $(\bbv^\star,\bbq^\star,\bblambda^\star)$ containing the optimal solution $(\bbv^\star,\bbq^\star)$ of \eqref{intvar} satisfies the following KKT conditions:
\begin{subequations}\label{eq:KKT_inc}\nonumber
	\begin{align}
	&\df(\bbv^\star)+\bbB^\T\bblambda^\star=\bb0, \label{eq:KKT_v_inc}\\
	&\dg(\bbq^\star)-\bblambda^\star+\partial\ind(\bbq^\star)\ni\bb0, \label{eq:KKT_q_inc}\\
	&\bbB\bbv^\star-\bbq^\star-\bbw=\bb0 \label{eq:KKT_lam_inc}
	\end{align}
\end{subequations}
where $\partial\ind(\bbq^\star)$ is the subdifferential set of $\ind(\cdot)$ at $\bbq^\star$. A subgradient of $\ind(\cdot)$ at any point $\bbq\in\ccalQ$, denoted as $\dind(\bbq)$, is defined as an element of the subdifferential set $\partial\ind(\bbq)$, i.e., $\dind(\bbq)\in\partial\ind(\bbq)$. Thus the above subdifferential inclusion condition \eqref{eq:KKT_inc} can be reinterpreted as that there exists $\dind(\bbq^\star)$ such that
\begin{subequations}\label{eq:KKT}
	\begin{align}
	&\df(\bbv^\star)+\bbB^\T\bblambda^\star=\bb0, \label{eq:KKT_v}\\
	&\dg(\bbq^\star)-\bblambda^\star+\dind(\bbq^\star)=\bb0, \label{eq:KKT_q}\\
	&\bbB\bbv^\star-\bbq^\star-\bbw=\bb0. \label{eq:KKT_lam}
	\end{align}
\end{subequations}	
In the sequel, we will always use this ``reinterpreted'' subgradient form to conduct the derivation since it is more straightforward. Such use of subgradient and notation have appeared in many recent references including \cite{bertsekas2011incremental,shi2015proximal}. When an equation involves a subgradient, which element in the corresponding subdifferential it uses is always clear from the context. For example, the subgradient used in \eqref{eq:KKT_q}, $\dind(\bbq^\star)$, equals to $\bblambda^\star-\dg(\tilde{\bbq}^\star)$ where
\[
\tilde{\bbq}^\star\triangleq\arg\min_{\bbq\in\ccalQ} g(\bbq)-\langle\bblambda^\star,\bbq\rangle.
\]	
Note that $\tilde{\bbq}^\star$ is well-defined when $g(\bbq)$ is strongly convex or $\ccalQ$ is compact, which is exactly the case in our paper. Under the aforementioned notations, our algorithm is equivalent to
\begin{equation}\label{eq:updates0}
\begin{array}{l}
\text{$\bbv$-update: }\bbv^{k+1}=\arg\min_\bbv f(\bbv)+ \langle\bblambda^k,\bbB\bbv\rangle;\\
\text{$\bbq$-update: }\bbq^{k+1}=\Proj{\bbq^k-\alpha\dg(\bbq^k)+\alpha\bblambda^k};\\
\text{$\bblambda$-update: }\bblambda^{k+1}=\bblambda^k+\beta(\bbB\bbv^{k+1}-\bbq^{k+1}-\bbw).\\
\end{array}
\end{equation}
The recursive relationship of the sequence $\{\bbv^k,\bbq^k,\bblambda^k\}$ in \eqref{eq:updates0} is further equivalent to\footnote{By using the notion of subgradient and indicator function, one is able to show that a gradient projection step can be understood as a subgradient step.}
\begin{subequations}\label{eq:updates1}
	\begin{align}
	&\df(\bbv^{k+1})+\bbB^\T\bblambda^k=\bb0,\label{eq:updates1_v}\\
	&\bbq^{k+1}=\bbq^k-\alpha\dg(\bbq^k)+\alpha\bblambda^k-\alpha\dind(\bbq^{k+1}),\label{eq:updates1_q}\\
	&\bblambda^{k+1}=\bblambda^k+\beta(\bbB\bbv^{k+1}-\bbq^{k+1}-\bbw).\label{eq:updates1_lam}
	\end{align}
\end{subequations}
The goal of convergence analysis is to show that the first-order residuals $\|\df(\bbv^{k+1})+\bbB^\T\bblambda^k\|$, $\|\dg(\bbq^k)-\bblambda^k+\dind(\bbq^{k+1})\|$, and $\|\bbB\bbv^{k+1}-\bbq^{k+1}-\bbw\|$ will all go to zero when $k$ goes to infinity. Due to the fact that the optimal solution exists and is unique, vanishing first-order residuals automatically imply the convergence of $(\bbv^k,\bbq^k)$ to $(\bbv^\star,\bbq^\star)$. Keeping this in mind, in the sequel, we will provide a theorem regarding convergence under suitable step-sizes.

\begin{theorem}
	\label{theo:con}
	Let $\seta$ and $\sL$ be the smallest and largest singular values of the reduced graph Laplacian matrix $\bbB$, respectively. If the non-negative step-sizes $\alpha$ and $\beta$ are chosen such that
	\begin{equation}\label{eq:step-size_rule}
	\left\{
	\begin{array}{lcl}
	\alpha&<&2/[\gamma(\sL^{-1}+\seta^{-1})], \\
	\beta&<&2/[\sL^{2}+\gamma^{-1}(\sL+\seta)],
	\end{array}
	\right.
	\end{equation}
	then the sequence $\{\bbq^k\}$ generated by the PPD updates (S1-S3) converges to the optimizer $\bbq^\star$.
\end{theorem}
\begin{IEEEproof}
	We first focus on showing that the successive difference $\{\|\bbq^k-\bbq^{k+1}\|_2^2$+$\|\bblambda^k-\bblambda^{k+1}\|_2^2\}$ is an infinitely summable sequence and thus converges to zero. Next, we prove that such vanishing successive difference implies vanishing first-order residuals. Therefore, it completes the proof.
	
	The first step would be connecting the $(k+1)$-th iterates to the corresponding optimal quantities. To this end, substituting \eqref{eq:updates1_lam} into \eqref{eq:updates1_v} and \eqref{eq:updates1_q}, respectively, for $\bblambda^k$ and then subtracting it from the KKT equalities \eqref{eq:KKT}, we have
	\begin{subequations}\label{eq:updates2}
		\begin{align}
		&\df(\bbv^{k+1})-\df(\bbv^\star)+\bbB^\T(\bblambda^{k+1}-\bblambda^\star) \nonumber \\
		&\qquad-\beta\bbB^\T(\bbB(\bbv^{k+1}-\bbv^\star)-(\bbq^{k+1}-\bbq^\star))=\bb0,\label{eq:updates2_v}\\
		&\bbq^{k+1}=\bbq^k-\alpha(\dg(\bbq^k)-\dg(\bbq^\star))+\alpha(\bblambda^{k+1}-\bblambda^\star)\nonumber\\
		&\qquad\quad-\alpha\beta(\bbB(\bbv^{k+1}-\bbv^\star)-(\bbq^{k+1}-\bbq^\star)) \nonumber \\
		&\qquad\quad-\alpha(\dind(\bbq^{k+1})-\dind(\bbq^\star)),\label{eq:updates2_q}\\
		&\bblambda^{k+1}=\bblambda^k+\beta(\bbB(\bbv^{k+1}-\bbv^\star)-(\bbq^{k+1}-\bbq^\star)).\label{eq:updates2_lam}
		\end{align}
	\end{subequations}
	Additionally, by the strong convexity and gradient Lipschitz continuity of $g(\cdot)$, we have
	\begin{equation}\label{eq:strong_and_Lip}
	\hspace{-1.5em}\begin{array}{rcl}
	&    &\frac{2\eta L}{\eta+L}\|\bbq^k-\bbq^\star\|_2^2+\frac{2}{\eta+L}\|\dg(\bbq^k)-\dg(\bbq^\star)\|_2^2\\
	&\leq&2\langle \bbq^k-\bbq^\star,\dg(\bbq^k)-\dg(\bbq^\star)\rangle\\
	&\leq&2\langle \bbq^{k+1}-\bbq^\star,\dg(\bbq^k)-\dg(\bbq^\star)\rangle\\
	&    &+\frac{\eta+L}{2}\|\bbq^k-\bbq^{k+1}\|_2^2+\frac{2}{\eta+L}\|\dg(\bbq^k)-\dg(\bbq^\star)\|_2^2.
	\end{array}
	\end{equation}
	Note that the first inequality of \eqref{eq:strong_and_Lip} is standard and obtained directly from the strong convexity and gradient Lipschitz continuity of $g(\cdot)$. We omit the proof as it is derived in Theorem 2.1.11 of \cite{nesterov2013introductory}.
	Reorganizing \eqref{eq:strong_and_Lip} gives
	\begin{equation}\label{eq:proof_2}
	\begin{array}{rcl}
	\frac{2\alpha\eta L}{\eta+L}\|\bbq^k-\bbq^\star\|_2^2
	&\leq& 2\langle \bbq^{k+1}-\bbq^\star,\alpha(\dg(\bbq^k)-\dg(\bbq^\star))\rangle\\
	&&+\frac{\alpha(\eta+L)}{2}\|\bbq^k-\bbq^{k+1}\|_2^2.
	\end{array}
	\end{equation}
	Substituting \eqref{eq:updates2_q} into \eqref{eq:proof_2} for $\alpha(\dg(\bbq^k)-\dg(\bbq^\star))$ leads to
	\begin{equation} %
	\hspace{-1.5em}\begin{array}{rcl}
	&    &\frac{2\alpha\eta L}{\eta+L}\|\bbq^k-\bbq^\star\|_2^2\\
	&\leq&2\langle \bbq^{k+1}-\bbq^\star,\bbq^k-\bbq^{k+1}\rangle+2\langle \bbq^{k+1}-\bbq^\star,\alpha\bblambda^{k+1}-\alpha\bblambda^\star\rangle\\
	&    &-2\langle \bbq^{k+1}-\bbq^\star,\alpha\beta \bbB(\bbv^{k+1}-\bbv^\star)\rangle\\
	&    &+2\langle \bbq^{k+1}-\bbq^\star,\alpha\beta(\bbq^{k+1}-\bbq^\star)\rangle\\
	&    &-2\langle \bbq^{k+1}-\bbq^\star,\alpha(\dind(\bbq^{k+1})-\dind(\bbq^\star))\rangle\\
	&    &+\frac{\alpha(\eta+L)}{2}\|\bbq^k-\bbq^{k+1}\|_2^2\\
	&\leq&2\langle \bbq^{k+1}-\bbq^\star,\bbq^k-\bbq^{k+1}\rangle+2\alpha\langle \bbq^{k+1}-\bbq^\star,\bblambda^{k+1}-\bblambda^\star\rangle\\
	&    &-2\alpha\beta\langle \bbq^{k+1}-\bbq^\star, \bbB(\bbv^{k+1}-\bbv^\star)\rangle+2\alpha\beta\|\bbq^{k+1}-\bbq^\star\|_2^2\\
	&    &+\frac{\alpha(\eta+L)}{2}\|\bbq^k-\bbq^{k+1}\|_2^2.
	\end{array} \label{eq:proof_3}
	\end{equation}
	Next, by the strong convexity of the function $f(\cdot)$, we have
	\begin{equation}\label{eq:proof_4}
	\begin{array}{c}
	2c\|\bbv^{k+1}-\bbv^\star\|_2^2\leq 2\langle \bbv^{k+1}-\bbv^\star, \df(\bbv^{k+1})-\df(\bbv^\star)\rangle.
	\end{array}
	\end{equation}
	Substituting \eqref{eq:updates2_v} into \eqref{eq:proof_4} for $\df(\bbv^{k+1})-\df(\bbv^\star)$ leads to
	\begin{equation}\label{eq:proof_5}
	\hspace{-1.5em}\begin{array}{rcl}
	&    &2\alpha c\|\bbv^{k+1}-\bbv^\star\|_2^2\\
	&\leq&2\alpha\langle \bbv^{k+1}-\bbv^\star, \bbB^\T(\bblambda^\star-\bblambda^{k+1})\\
	&    &+\beta \bbB^\T(\bbB(\bbv^{k+1}-\bbv^\star)-(\bbq^{k+1}-\bbq^\star))\rangle\\
	&  = &2\alpha\langle \bbB(\bbv^{k+1}-\bbv^\star), \bblambda^\star-\bblambda^{k+1}\rangle+2\alpha\beta\|\bbv^{k+1}-\bbv^\star\|_{\bbB^\T \bbB}^2 \\
	&    &-2\alpha\beta\langle \bbB(\bbv^{k+1}-\bbv^\star),\bbq^{k+1}-\bbq^\star\rangle
	\end{array}
	\end{equation}
	where the weighted norm $\|\bby\|_\bbB ^2 := \bby^\T\bbB\bby$ for any vector $\bby$. Summing up \eqref{eq:proof_3} and \eqref{eq:proof_5} results in
	\begin{equation}\label{eq:proof_6}
	\hspace{-1.5em}\begin{array}{rcl}
	&    &\frac{2\alpha\eta L}{\eta+L}\|\bbq^k-\bbq^\star\|_2^2+2\alpha c\|\bbv^{k+1}-\bbv^\star\|_2^2\\
	&\leq&2\langle \bbq^{k+1}-\bbq^\star,\bbq^k-\bbq^{k+1}\rangle+2\alpha\langle \bbq^{k+1}-\bbq^\star,\bblambda^{k+1}-\bblambda^\star\rangle\\
	&    &-2\alpha\beta\langle \bbq^{k+1}-\bbq^\star, \bbB(\bbv^{k+1}-\bbv^\star)\rangle+2\alpha\beta\|\bbq^{k+1}-\bbq^\star\|_2^2\\
	&    &+\frac{\alpha(\eta+L)}{2}\|\bbq^k-\bbq^{k+1}\|_2^2\\
	&    &+2\alpha\langle \bbB(\bbv^{k+1}-\bbv^\star), \bblambda^\star-\bblambda^{k+1}\rangle+2\alpha\beta \|\bbv^{k+1}-\bbv^\star\|_{\bbB^\T \bbB}^2\\
	&    &-2\alpha\beta\langle \bbB(\bbv^{k+1}-\bbv^\star),\bbq^{k+1}-\bbq^\star\rangle.
	\end{array}
	\end{equation}
	Applying the basic inequality 
	\[
	2\langle \sqrt{\rho}\bba,\sqrt{\rho^{-1}}\bbb\rangle\leq \rho\|\bba\|^2+\rho^{-1}\|\bbb\|^2,
	\]
	which holds for any $\rho>0$ and any real vectors $\bba$ and $\bbb$ of the same dimension and utilizing \eqref{eq:updates2_lam}, the right-hand-side of \eqref{eq:proof_6} can be relaxed as
	\begin{equation}\label{eq:proof_6_right}
	\hspace{-1.5em}\begin{array}{rl}
	&2\langle \bbq^{k+1}-\bbq^\star,\bbq^k-\bbq^{k+1}\rangle+\frac{2\alpha}{\beta}\langle \bblambda^{k}-\bblambda^{k+1},\bblambda^{k+1}-\bblambda^\star\rangle\\
	&+(\frac{(1-\varepsilon)\alpha}{\beta}+\frac{(1+\varepsilon)\alpha}{\beta})\|\bblambda^k-\bblambda^{k+1}\|_2^2+\frac{\alpha(\eta+L)}{2}\|\bbq^k-\bbq^{k+1}\|_2^2 \\
	&\leq\|\bbq^k-\bbq^\star\|_2^2-\|\bbq^{k+1}-\bbq^\star\|_2^2-\|\bbq^k-\bbq^{k+1}\|_2^2\\ &+\frac{\alpha}{\beta}(\|\bblambda^{k}-\bblambda^\star\|_2^2-\|\bblambda^{k+1}-\bblambda^\star\|_2^2)-\varepsilon\|\bblambda^{k}-\bblambda^{k+1}\|_2^2\\
	&+\frac{\alpha(\eta+L)}{2}\|\bbq^k-\bbq^{k+1}\|_2^2+(1+\varepsilon)\alpha\beta(1+\rho)\|\bbB(\bbv^{k+1}-\bbv^\star)\|_2^2\\
	&+(1+\varepsilon)\alpha\beta(1+\frac{1}{\rho})\|\bbq^{k+1}-\bbq^\star\|_2^2
	\end{array}
	\end{equation}
	where $\varepsilon$ is any arbitrary constant in the interval $(0,1)$. It then follows from \eqref{eq:proof_6} and \eqref{eq:proof_6_right} that
	\begin{equation}\label{eq:proof_7}
	\hspace{-1em}\begin{array}{rl}
	&(1-\frac{\alpha(\eta+L)}{2})\|\bbq^k-\bbq^{k+1}\|_2^2+\varepsilon\|\bblambda^{k}-\bblambda^{k+1}\|_2^2\\
	\leq&(1-\frac{2\alpha\eta L}{\eta+L})\|\bbq^k-\bbq^\star\|_2^2+\frac{\alpha}{\beta}(\|\bblambda^{k}-\bblambda^\star\|_2^2-\|\bblambda^{k+1}-\bblambda^\star\|_2^2)\\
	&-(1-(1+\varepsilon)\alpha\beta(1+\frac{1}{\rho}))\|\bbq^{k+1}-\bbq^\star\|_2^2\\
	&-(2\alpha c-(1+\varepsilon)\alpha\beta(1+\rho)\sigmax{\bbB^\T \bbB})\|\bbv^{k+1}-\bbv^\star\|_2^2
	\end{array}
	\end{equation}
	where $\sigmax{\cdot}$ measures the largest singular value of a matrix. In order to conclude that $\{\|\bbq^k-\bbq^{k+1}\|_2^2+\|\bblambda^k-\bblambda^{k+1}\|_2^2\}$ is infinitely summable, based on \eqref{eq:proof_7}, the following conditions must be met (this is done by comparing the coefficients of terms in \eqref{eq:proof_7} and requiring non-positive coefficients on the right-hand-side of \eqref{eq:proof_7} after telescope cancellation):
	\begin{equation}\label{eq:proof_8}
	\left\{
	\begin{array}{l}
	1-\frac{\alpha(\eta+L)}{2}>0, \\
	\varepsilon>0, \\
	1-\frac{2\alpha\eta L}{\eta+L}\leq1-(1+\varepsilon)\alpha\beta(1+\frac{1}{\rho}), \\
	2\alpha c-(1+\varepsilon)\alpha\beta(1+\rho)\sigmax{\bbB^\T \bbB}\geq0, \\
	\end{array}
	\right.
	\end{equation}
	which is equivalent to requiring
	\begin{equation}\label{eq:proof_9}
	\left\{
	\begin{array}{cl}
	\alpha&<\frac{2}{\eta+L}, \\
	\beta&<\min\left\{\frac{2\eta L}{(\eta+L)(1+\frac{1}{\rho})},\frac{2c}{(1+\rho)\sigmax{\bbB^\T \bbB}} \right\}.\\
	\end{array}
	\right.
	\end{equation}
	To maximize the possible range of $\beta$, we choose $\rho=\frac{c(\eta+L)}{\eta L\sigmax{\bbB^\T \bbB}}$. Accordingly, we obtain the step-size rule
	\begin{equation}\label{eq:proof_10}
	\left\{
	\begin{array}{cl}
	\alpha&<\frac{2}{\eta+L}, \\
	\beta&<\frac{2c\eta L}{\eta L \sigmax{\bbB^\T \bbB}+c(\eta+L)}.\\
	\end{array}
	\right.
	\end{equation}

	Under the choice of \eqref{eq:proof_10}, we can find the infinite summability of $\{\|\bbq^k-\bbq^{k+1}\|_2^2+\|\bblambda^k-\bblambda^{k+1}\|_2^2\}$ from \eqref{eq:proof_7}. Hence, $\lim_{k\rightarrow\infty}\|\bbq^k-\bbq^{k+1}\|_2=0$ and $\lim_{k\rightarrow\infty}\|\bblambda^k-\bblambda^{k+1}\|_2=0$. Noticing the relation \eqref{eq:updates1}, we thus have
	\begin{subequations}\label{eq:limit}
		\begin{align}
		&\|\df(\bbv^{k+1})+\bbB^\T\bblambda^k\|_2=0, \label{eq:limit_v}\\
		&\lim_{k\rightarrow\infty}\|\dg(\bbq^k)-\bblambda^k+\dind(\bbq^{k+1})\|_2=0, \label{eq:limit_q}\\
		&\lim_{k\rightarrow\infty}\|\bbB\bbv^{k+1}-\bbq^{k+1}-\bbw\|_2=0. \label{eq:limit_lam}
		\end{align}
	\end{subequations}
	Comparing \eqref{eq:limit} and \eqref{eq:KKT}, $(\bbv^\infty,\bbq^\infty,\bblambda^\infty)$ satisfies the KKT conditions \eqref{eq:KKT}. Therefore, we conclude $\{\bbq^k\}$ converges to $\bbq^\star$.
	
	The explicit forms of $g(\cdot)$ and $f(\cdot)$ are $g(\bbq)=\frac{\gamma}{2} \|\bbX (\bbq+\bbw)-\bbmu\|_\bbB^2$ and $f(\bbv)=\frac{1}{2} \|\bbv - \bbmu\|^2$. The Hessian of $g(\cdot)$ is $\nabla^2g=\gamma\bbX$ while that of $f(\cdot)$ is an identity matrix. Let us denote $\seta$ and $\sL$ as the smallest and largest singular values of the reduced Laplacian matrix $\bbB$, respectively. Noticing that $\bbX=\bbB^{-1}$, we have $\eta=\gamma\sL^{-1}$ and $L=\gamma\seta^{-1}$. Obviously $c=1$ in this case. Finally, substituting the values of $\eta$, $L$, and $c$ into \eqref{eq:proof_10} leads to the specific bounds stated in the theorem, i.e., \eqref{eq:step-size_rule}.
\end{IEEEproof}

{\ml It is clear that $\alpha=O(\seta)$ and $\beta=O(1/\sL^2)$. Increasing the network size which leads to the growing of $\sL$ may reduce the bound on $\beta$. Nonetheless, the dual step-size $\beta$ does not significantly affect the practical convergence speed empirically.  In contrast, we usually find the primal step-size $\alpha$ playing a more crucial role in adjusting the convergence speed. In addition, the quantity $\seta$ ``sort of''\footnote{We comment it as ``sort of'' because conventionally, the algebraic connectivity of the graph is defined as the second smallest eigenvalue of the standard graph Laplacian. There is some relationship between the algebraic connectivity and $\seta$, but it is beyond the scope of our discussion.} entails the connectivity of the network since it is the smallest eigenvalue of the reduced graph Laplacian. But as long as the network is connected, this quantity is lower bounded away from zero. Empirically, a larger step-size bound implies the possibility of a faster convergence. Accordingly, a better connected physical network would result in a faster convergence, which aligns with our common intuition.}
	
\begin{remark} \label{rmk:hubber} (Relaxations and generalizations) 
\ws{Our convergence analysis is conducted under the assumption of functions $f$ and $g$ both being strongly convex. This assumption can be relaxed to the so-called restricted strong convexity (Huber loss satisfies such assumption) \cite{shi2015extra}. In this case, the first inequality in \eqref{eq:strong_and_Lip} needs to be replaced by a looser one (worse coefficients) and a narrower stable step-size region will be derived.}

\ws{In addition, it is possible to totally remove the (restricted) strong convexity assumption on $f$ with a simple modification over our current scheme. Specifically, the \eqref{v-update} needs to be tweaked by appending a proximal term, as given by}
\[
\bbv^{k+1}=\arg\min_\bbv f_1(\bbv)+ \langle\bblambda^k,\bbB\bbv\rangle+\vartheta\|\bbv-\bbv^{k}\|^2.
\]
\ws{Informally speaking, adding such proximal term $\vartheta\|\bbv-\bbv^{k}\|^2$ enhances the stability of the algorithm when dealing with non-strongly convex functions. The final step-size rule for $\alpha$ and $\beta$ would foreseeably not be affected by the strong convexity constant $c$, but rather depend on $\vartheta$. Consequently, $\vartheta$ needs to be adjusted together with $\alpha$ and $\beta$ in the PPD algorithm to satisfy some specific rule. However, adding a proximal term could presumably slow down the convergence speed.}
\end{remark}

\subsection{Online Feedback Design}
Thus far, we assume the availability of static $w_j$ per bus $j$. However, it would change in accordance with the system operating condition. To account for system dynamics based on \eqref{ldfVm2}, we need to compute and update $w_j$ which requires network-wide complex power injections. To this end, two-way communications between each bus and a centralized computer are necessary. Generally, this is not feasible due to limiting communication in distribution network in addition to fast system dynamics. Thanks to the sparsity of $\bbB$, its entry is zero for any pair of buses whose corresponding buses are not connected by a line segment. Accordingly, we propose to obtain a time-varying $w_j^k$ through neighboring voltage measurement exchanges. Each bus $j$ measures its voltage magnitude $\tv_j^{k+1}$ after (S2) and broadcast to its neighboring buses. Knowing incident line reactance values at each bus, we update $w_j^{k+1}$ as
%{\ml We assume that inverter controllers have some basic functionality of computing simple updates and communicating infrequently with a centralized supervising computer. Each bus would be able to obtain the incident line reactance values from a centralized supervising computer, e.g., distribution management system (DMS). Therefore, if the network topology varies due to some power system operations, the supervising computer would be able to construct a new matrix $\bbB$ corresponding to system changes and send necessary line reactance values to each controller. Knowing incident line reactance values at each bus}, we update $w_j^{k+1}$ as
%
\begin{equation}
\label{wupdate}
w_j^{k+1}:=\sum_{i \in \mathcal{N}_{j}} B_{ji} \tv_i^{k+1} - q_j^{k+1}
\end{equation}
which is computed locally and adopted in (S3). Note that each bus would only need to store the neighboring line reactance values to update this feedback signal. Therefore, the memory requirement is minimal. The attractive features of the proposed feedback design are three-fold. First, we can obtain $w_j^{k+1}$ locally at each bus by adopting bus-to-bus communication architecture. This is nicely designed to our HVC scheme as the distributed feature is maintained. Second, the instantaneous voltage measurements contain the latest system information, and hence $\bbw^{k+1}$ effectively approximates the dynamically varying operating conditions. Last but not least, the voltage feedback control design improves the robustness to mismatch and imperfection in system modeling (see e.g., \cite[Sec. 8.9]{luenberger_dyn79}) since the voltage measurements could potentially capture the underlying non-linearity in the power networks.
{\ml
\begin{remark} \label{rmk:nodetonode} (Cyber Network Topology) {\rm
	The node-to-node architecture of the HVC design can be generalized to instead coordinate clusters of buses as long as the cyber network is connected; see e.g., \cite{dall_tsg13}. This way, it is not necessary to have DERs to be connected to each other by a line segment. Additionally, even if the distribution network is not a complete entity (i.e., a DER is not necessarily attached to every bus), we may eliminate all the buses with no DERs to create an equivalent network by adopting the Kron Reduction method \cite{ARBergen2000}. Thus, this reduced network consists of only buses with DERs installed. By adopting the voltage-based feedback signal $\bbw$ in \eqref{wupdate}, the corresponding reduced {\it Bbus matrix} $\bbB$ would explicitly account for all system characteristics as to the original network. Accordingly, the performance of our HVC design can still be guaranteed under a generalized distribution network where DERs are not attached to each and every bus.
}
\end{remark}}
\subsection{Limited Communication Rates}
\label{sec:AHVC}
The performance of the proposed HVC design relies on the quality of bus-to-bus communication links, which we have assumed to be perfect throughout the algorithmic design. However, random link failures and messaging delays are common because of either network congestion, or poor signal-to-noise ratios in some wireless environments for a contemporary digital communication system. It is imperative to examine how the PPD-based HVC scheme works under imperfect communication, which leads to the following two different scenarios. One is often referred to as asynchronous networking that is consisting of both link failures and messaging delays \cite{chang2015asynchronous,peng2015arock}. Meanwhile, the other only considers link failures. The later one can also be referred to as time-varying network \cite{hong2016stochastic,nedich2016achieving}. Albeit the second scenario seems to be a special case of the first one, by playing a simple trick of embedding a time stamp in the message between each pair of buses, delays can also be treated as link failures\footnote{We assume that each bus has a clock that is aligned. As the information exchanges and dual updates could be performed at a relatively low speed, a slight mismatch between clocks does not break the viability of this approach.}. Informally speaking, as long as the delay is bounded, and the time-varying communication network is $B$-connected\footnote{It is a connectivity description of graphs under time-varying scenarios. Readers are referred to Assumption 2 at page 7 of reference \cite{nedich2016achieving} for detailed definition.}, one should be able to choose small enough step-sizes to stabilize the proposed algorithm. There have been some analysis for algorithms under these conditions in the literature (see, e.g., \cite{iutzeler2013asynchronous,bianchi2014mlsp}). Rigorous proof of convergence properties under these conditions is out of the scope of this paper and will be a future direction. We have tested and validated the proposed design using realistic distribution networks in Sec. \ref{sec:num}.

To tackle the challenge posed by imperfect communication networks, we leverage the work in \cite{iutzeler2013asynchronous,bianchi2014mlsp,Zhu_tsp09} regarding the "freezing" strategy for distributed optimization problems. Conceptually, every PPD variable remains unchanged until new information becomes available from neighboring buses. The asynchronous version of a related distributed primal-dual algorithm in \cite{bianchi2014mlsp} has been proven to be convergent under random activation of agents, i.e., link $(i,j)\in\ccalE$ is available only when $i$ and $j$ are both randomly activated. It is assumed that the activation of each bus follows a Bernoulli distribution, independent across time. Nonetheless, under this strategy, the aforementioned PPD-based HVC updates would completely halt under the case of a total link failure in the communication network. The novelty of our work lies in the extension of the HVC scheme by modifying update steps (S1)-(S3) to have: a) satisfactory performance under partial link failures, b) capability to continue providing VAR regulation under a total link failure scenario. To this end, we freeze the variables $v_j \textrm{ and } \lambda_j$ associated with the inactive bus $j$ while always updating the VAR control signal according to \eqref{qjupdate} by adapting local gradient information from $\tv_j^k$. Since the voltage measurements always contain the most updated information of the network, the local voltage control design objective $f_2$ is advocated to continue providing the VAR support. As a result, under a total communication failure, our HVC framework boils down to a surrogate local controller design based on the current value of $\lambda_j$. This is similar to the microgrid secondary frequency/voltage control design where $\lambda_j$ can be treated as an offset signal to a local droop controller \cite{lu_naps16,lu_tetcas15}. We denote $\ccalN_a^k \subseteq \ccalN$ as the subset activated nodes at iteration $k$. Our proposed asynchronous (A-)HVC algorithm is tabulated in Algorithm \ref{alg:linkfailure}. 
%The attractive feature of this asynchronous (A-)HVC architecture are three-fold. First, the algorithm does not need to wait for convergence to determine the control input $\bbq^k$. 
%
\begin{algorithm}[t]
	\caption{Asynchronous HVC (A-HVC) algorithm} %under communication link failure}
	\label{alg:linkfailure}
	\begin{algorithmic}[1]
		\For{every iteration $k = 1, 2, \ldots$}
		\For{bus $j \in \ccalN_a^k$}
		\State (AS1): update $v_j^{k+1}$ as in \eqref{vbarup};
		\State (AS2): update $q_j^{k+1}$ as in \eqref{qjupdate};
		\State Update $w_j^{k+1}$ as in \eqref{wupdate};
		\State (AS3): update $\lambda_j^{k+1}$ as in \eqref{lambup};
		\EndFor
		\For{bus $j \notin \ccalN_a^k$}
		\State $v_j^{k+1}=v_j^{k}$;
		\State (AS2): update $q_j^{k+1}$ as in \eqref{qjupdate};
		\State $w_j^{k+1}=w_j^{k}$
		\State $\lambda_j^{k+1}=\lambda_j^{k}$
		\EndFor
		\EndFor
	\end{algorithmic}
\end{algorithm}
%
%%%%%%%%%%%%%%%%%%%%%%%%%%%%%%%%%%%%%%%%%%%%%%%%%%%%%%%%%%%%%%%%%%%%%
%%
%%      Section: Numberical Tests %%
%%%%%%%%%%%%%%%%%%%%%%%%%%%%%%%%%%%%%%%%%%%%%%%%%%%%%%%%%%%%%%%%%%%%
%%
\section{Numerical Tests} \label{sec:num}
The numerical tests presented in this section demonstrate the effectiveness of the proposed communication-cognizant HVC design for practical distribution feeders. We investigate the performance of our scheme under the settings of both static and dynamically time-varying network operating conditions. The dynamic tests are performed using the IEEE 123-bus test case \cite{ieee123}. The desired voltage magnitude $\mu_j$ is chosen to be $1$ p.u. at every bus $j$. Furthermore, each bus is assumed to have a certain number of PV panels installed, and thus it is able to control/provide VAR via advanced inverter design. Albeit the HVC design is based on the linearized model \eqref{ldfVm2}, we test and validate the performance using the full ac power flow model. All numerical tests are performed using MathWorks\textsuperscript{\textregistered} MATLAB 2014a software and the OpenDSS for solving the actual power flow. Accordingly, the bus voltage magnitude {\ml $\tdbbv$}, instead of the one obtained from \eqref{ldfVm2}, is used for VAR control outputs in \eqref{qjupdate} and the following numerical tests. 

\subsection{Static System Operating Conditions}
\begin{figure}[tb]
	\centering
	\includegraphics[width=0.6\linewidth,clip = true, trim = 0.3in 0.3in  .5in  0.16in]{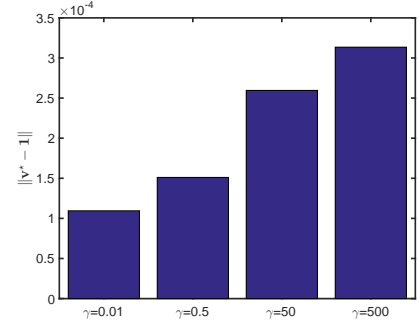}
	\caption{Voltage mismatch norm error versus various values of importance factor $\gamma$ across the network under a static system setting.}
	\label{fig:gamma_bar}
	\vspace{-4mm}
\end{figure}
{\ml A single-phase radial power distribution feeder that consists of 21 buses with $v_0=1$ at the head of the feeder is first adopted to test the algorithm under the static setup. %This network %is equivalent to the system in Fig. \ref{fig:radial} for has $N=20$ with t
The impedance of each line segment is set to be  $(0.233+j0.366)\Omega$.} Per bus $j$, we fix the loading $p^\mathrm{c}_j=70\textrm{kW} \textrm{ and } q^\mathrm{c}_j=20\textrm{kVAR}$ while choosing the inverter rating to be  $(70+\psi)$kVA where $\psi$ is zero-mean Gaussian having variance  13.33, thus modeling the variation in inverter sizing by $50\%$. Accordingly, the VAR constraints in \eqref{invlim} would become active at some locations. We test the HVC algorithm with various choices of importance factor $\gamma$. To demonstrate the trade-off between distributed and local control designs, we plot in Fig. \ref{fig:gamma_bar} the optimal voltage mismatch norm error $\|\bbv^\star-\mathbf{1}\|$ for each $\gamma$ value with all other settings the same.  Notice that increasing $\gamma$ value in term of adding more weight on the local control objective $f_2$ results in a larger voltage mismatch error. This corroborates with our earlier claim that local control schemes attain a sub-optimal VAR setting under limited VAR resources. Given a system model, one may study this trade-off offline to tune the importance factor $\gamma$ accordingly. 
%For this feeder, we let $\gamma=0.5$.
%
\begin{figure}[tb]
	\centering
	\includegraphics[width=0.8\linewidth,clip = true, trim = 0.5in .2in  .6in  .1in]{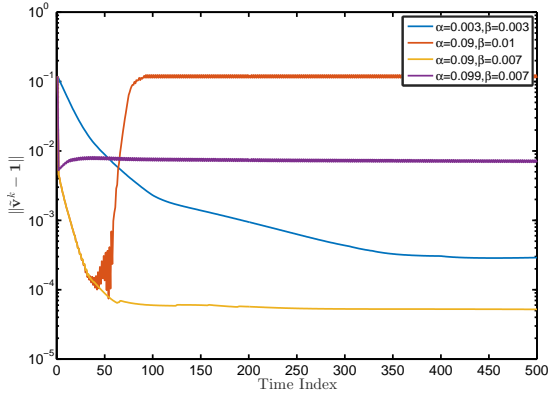}
	\caption{Log-scale voltage mismatch norm error versus the total number of updates across the network with different step-size choices of $\alpha \textrm{ and } \beta$ under the static system setting.}
	\label{fig:stepsize}
	\vspace{-4mm}
\end{figure}

Based on the convergence properties in Theorem \ref{theo:con} with $\gamma=0.5$, we have $\alpha<0.092 \textrm{ and } \beta<0.0073$. Fig. \ref{fig:stepsize} plots the iterative voltage mismatch norm error {\ml $\|\tdbbv^k-\mathbf{1}\|$} in log-scale for various step-size choices assuming a perfect communication. To violate the steps-size constraints, we let $\alpha=0.099$ and $\beta=0.01$ in two different scenarios, respectively. It clearly shows that the HVC design fails to converge under these cases. To {\ml stabilize} our design, we bound the step-size values to be within their limits as depicted in Fig. \ref{fig:stepsize}. Note that the effect of step-size choices shows a trade-off between the stability and convergence rate. The larger $\alpha \textrm{ and } \beta$ are, the faster the updates converges. Nonetheless, this could potentially lead to oscillations in the error performance, exhibiting instability under fast dynamics. To tackle this problem, once the full feeder information becomes available, our convergence properties in Theorem \ref{theo:con} are very useful in terms of selecting proper step-size choices. Otherwise, it is also possible to adjust the step-size on-the-fly by decreasing the values based on its local voltage oscillation intensity. To sum up, under appropriate step-size choices, Fig. \ref{fig:stepsize} validates the effectiveness of our scheme, in terms of achieving the optimal VAR setting while requiring no centralized coordination.

\begin{figure}[tb]
	\centering
	\includegraphics[width=0.8\linewidth,clip = true, trim =  0.5in .26in  .6in  .1in]{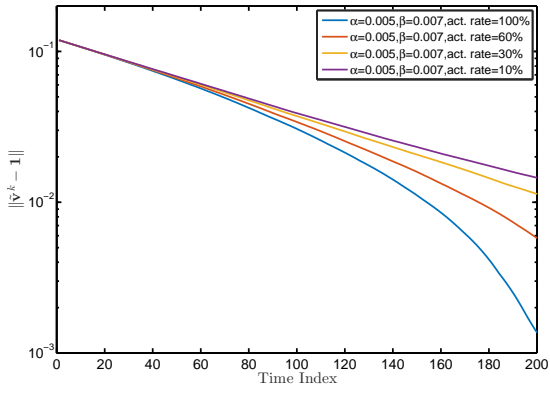}
	\caption{Log-scale voltage mismatch norm error versus the total number of updates across the network with fixed step-size choices of $\alpha \textrm{ and } \beta$ and varying bus activation rate under the static system setting.}
	\label{fig:act_rate}
	\vspace{-4mm}
\end{figure}
Moreover, to showcase the robustness of our A-HVC scheme under imperfect communication links, we model the activation of every bus $j$ as a Bernoulli distribution with the same probability to each other. Since the distributed part of the HVC are most likely to be affected by random communication link failures, we let $\beta > \alpha$ to investigate the performance of A-HVC design. To this end,  we have $\alpha=0.005, \beta=0.007, \textrm{ and } \gamma=0.05$ while fixing other settings to be the same as the earlier test. Fig. \ref{fig:act_rate} plots the iterate voltage mismatch error in log-scale under various bus activation rate ranging from $10\% \textrm{ to } 100\%$, where the case of $100\%$ corresponds to the perfect communication scenario (synchronous case). It clearly depicts that our design enjoys a satisfactory performance guarantee under random link failures for regulating the network voltage. Informally speaking, a lower bus activation rate would lead to a slower convergence speed, with a no link failure scenario exhibiting the fastest convergence. This test verifies that our proposed A-HVC design is robust against imperfect communication and thus able to cope with cyber resource constraints.
\begin{figure}[tb]
	\centering
	\includegraphics[width=0.8\linewidth, height = 3.8cm, clip = true, trim = 0.6in 0.1in  .6in  0.1in]{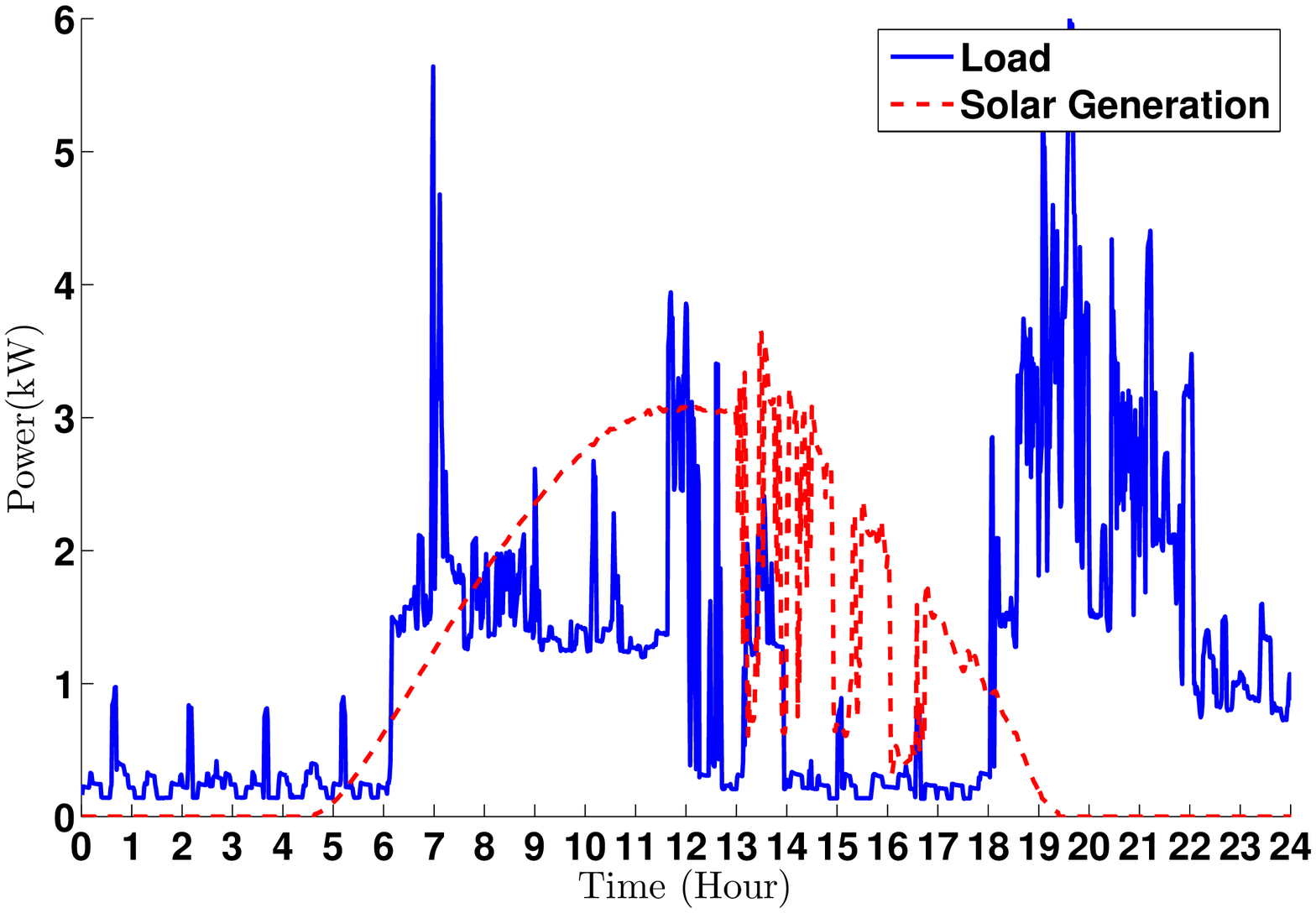}
	\caption{Sample daily load and solar generation profiles of a single residential home.}
	\label{fig:load_profile}
	\vspace{-4mm}
\end{figure}

\subsection{Dynamic System Operating Conditions}
To corroborate our HVC scheme for online implementation, we have tested the proposed algorithm on the IEEE 123-bus test case \cite{ieee123}. Dynamic system operating conditions are generated using real load profiles from an online data repository \cite{UCI_data}, as shown in Fig. \ref{fig:load_profile}. The minute-sampled active and reactive power consumption data along with solar profile were collected on Friday, June 20, 2010  for residential loads. Each home shares similar load and generation patterns, which are diversified by small random additive noises. For each load node of the 123-bus test case, it is attached to a certain number of residential homes with solar generation rated at 3.5kW peak capacity. Additionally, for each minute time slot, physical VAR limits $[\underline{\bbq},\overline{\bbq}]$ are updated according to their inverter ratings (i.e., 3.5kVA per inverter) and instantaneous active power from solar generations. Under these settings, it turns out that VAR resources of inverters are insufficient to achieve perfectly flat voltage at all times. This implies the VAR constraints in problem \eqref{intvar} are active.

%{\ml We assume that each residential home shares the same load pattern depending on the hour of a day. Similarly, the PV generation pattern should be uniform across the network due to the proximity of home locations. Albeit each home follows similar load and generation patterns, these patterns are diversified by small random additive noises. For each load node of the 123-bus test case, it is attached to a certain number of residential homes with solar generation rated at 3.5kW peak capacity. The number of homes for each node is determined by rounding up the ratio between the active spot load demand in the original case information and the maximum daily active load (6kW) of a residential home (as depicted in Fig \ref{fig:load_profile}). This way, the location of DERs corresponds to nodes with active power loading in the IEEE 123-bus test case. For each minute time slot, physical VAR limits $[\underline{\bbq},\overline{\bbq}]$ are updated according to their inverter ratings (i.e., 3.5kVA per inverter at a home) and instantaneous active power from solar generations.} {\ml We also fix the tap positions of voltage regulators in order to better capture the performance of different inverter-based voltage control designs. Under these settings, it turns out VAR resources of inverters are insufficient to achieve perfectly flat voltage at all times. This implies the VAR constraints in problem \eqref{intvar} are active.}
%
\begin{figure}[tb]
	\centering
	\includegraphics[width=\linewidth,clip = true, trim =  0.5in .26in  .6in  .1in]{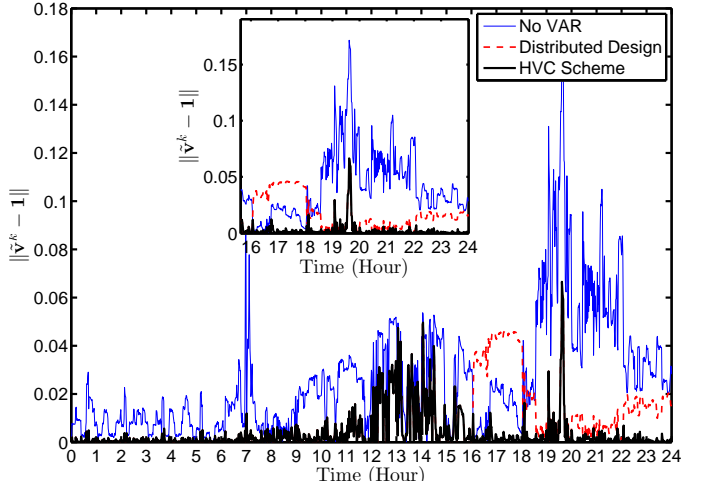}
	\caption{Daily voltage mismatch error for a-phase of the 123-bus under three different control strategies. A total communication link failure occurs from hour 16:00 to 24:00.}
	\label{fig:dynamic}
	\vspace{-4mm}
\end{figure}
{\ml Fig. \ref{fig:dynamic} plots the daily voltage mismatch for the a-phase of the 123-bus, where the other two phases exhibit a similar comparisons. Three different control strategies including no VAR support, distributed design, and our HVC scheme are plotted.} For the benchmark case of no VAR support, there are some under- and over-voltage issues due to load and solar variations. To improve the voltage quality, the proposed Algorithm \ref{alg:linkfailure} is tested for online implementation where every PPD-based iteration of (S1)-(S3) updates every 2 seconds assuming a fixed loading during corresponding one-minute interval. Additionally, we validate the robustness of the HVC design to the worst-case scenario of a total communication link outage from hour 16:00 to 24:00. {\ml Clearly, during the noon hours, that all three scenarios have similar voltage mismatch error because of limited VAR capability.} For all other hours, especially during the evening hours in the zoom-in view, we see that our HVC design effectively minimizes network voltage mismatch error, maintaining a nearly flat voltage profile. The attractive feature of our communication-cognizant HVC design enjoys a satisfactory performance even under the total link failure scenario, whereas the distributed-based control updates come to a complete halt, i.e. setting $q_j^{k+1}=q_j^k, \forall j \notin \ccalN_a^k$, resulting in a higher voltage mismatch error around hour 16:00-18:00 compared to the ones in the benchmark and HVC cases. 
{\ml To sum up, the proposed HVC design can efficiently regulate the voltage level by coordinating network-wide VAR resources. Meanwhile, its cognizant feature to the instantaneous availability of communication links is also attractive, considering the limited deployment of cyber infrastructure in distribution networks. Therefore, the proposed HVC design would facilitate the future engagements in inverter-based VAR resources to improve voltage support by accounting for practical constraints in both physical and cyber layers. }

%%%%%%%%%%%%%%%%%%%%%%%%%%%%%%%%%%%%%%%%%%%%%%%%%%%%%%%%%%%%%%%%%%%%%
%%
%%      Section: Conclusions %%
%%%%%%%%%%%%%%%%%%%%%%%%%%%%%%%%%%%%%%%%%%%%%%%%%%%%%%%%%%%%%%%%%%%%
%%
\section{Conclusions} \label{sec:con}
This paper has developed the communication-cognizant hybrid voltage control (HVC) scheme to coordinate network-wide VAR support in power distribution networks. {\ml Considering limited VAR resources, we have cast the specially designed voltage control problem combining both attractive features of distributed and local control architectures.} The PPD-based algorithm is evoked and only requires voltage measurement exchanges among neighboring buses with local computations. Moreover, we have provided the convergence properties of the aforementioned algorithm for proper step-size choices. To cope with cyber resource constraints and lack of reliable communication links, we have further extended the HVC design to have robustness against random communication link failures and, in particular, communication-cognizant feature to account for the worst-case scenario of a total communication outage. We have extensively validated the effectiveness of the HVC design using realistic distribution networks under both static and dynamic testing environments.

For future work, we plan to study the convergence properties of the proposed A-HVC design while offering the performance analysis of the online implementation under stochastic settings. We are also interested in other cyber-security aspects of inverter control designs.

%%%%%%%%%%%%%%%%%%%%%%%%%%%%%%%%%%%%%%%%%%%%%%%%%%%%%%%%%%%%%%%%%%%%%%
%                                                                    %
%       Bibliography %
%
%%%%%%%%%%%%%%%%%%%%%%%%%%%%%%%%%%%%%%%%%%%%%%%%%%%%%%%%%%%%%%%%%%%%%%
%\begin{thebibliography}{1}

%\end{thebibliography}

%%\newpage
%\section*{Biographies}
%
%\begin{IEEEbiography}[{\includegraphics[width=1in,
%height=1.25in,keepaspectratio]{haojanliu.eps}}]{Hao Jan (Max) Liu}
%(M'12) is currently a graduate student of ECE at the University of Illinois, Urbana-Champaign (UIUC). 
%He received her B.S. from Missouri Science and Technology, Rolla, MO in 2011, and the M.Sc. and Ph.D. degrees from the University of Illinois at Urbana Champaign, IL in 2013 and 2016, respectively. 
%His current research interests include power system monitoring and operations, dynamics and stability, distribution systems, and energy data analytics.
%\end{IEEEbiography}

\bibliographystyle{IEEEtran}
\bibliography{IEEEabrv,VAR,hzpub}

\end{document}